\documentclass[12pt]{amsart}
\usepackage{a4wide}
\usepackage{amssymb}
\usepackage{amsthm}
\usepackage{amsmath}
\usepackage{amscd}
\usepackage{verbatim}
\usepackage{color}
\usepackage[mathscr]{eucal}
\textheight8.8in \textwidth6.55in \numberwithin{equation}{section}

\newtheorem{theorem}{Theorem}

\newtheorem*{conj}{\bf Conjecture}

\theoremstyle{remark}

\numberwithin{equation}{section}

\newcommand{\Z}{\mathbb{Z}}

\begin{document}
\title{On a conjecture of Berndt and Kim}

\author{Kathrin Bringmann}
\address{Mathematical Institute,  University of Cologne, Weyertal 86-90, 50931 Cologne,
Germany} \email{kbringma@math.uni-koeln.de}
\author{Amanda Folsom}
\address{Yale University, Mathematics Department,
New Haven, CT 06520} \email{amanda.folsom@yale.edu}

\maketitle
\begin{abstract} We prove a recent conjecture of Berndt and Kim regarding the positivity of the coefficients in the asymptotic expansion of a class of partial theta functions.  This generalizes results found in RamanujanÕs second notebook, and recent work of Galway and Stanley.
\end{abstract}
\section{Introduction}\label{intro}

In his second notebook \cite{Be}, page 324, Ramanujan claims an asymptotic expansion for the partial theta function
\begin{equation}\label{Ram}
2\sum_{n=0}^\infty (-1)^n q^{n^2+n}\sim 1+t+t^2+2t^3+5t^4+\ldots
\end{equation}
with $q=\frac{1-t}{1+t}$ as $t\to 0^+$. It is not clear, just given the left-hand side of \eqref{Ram}, that the coefficients of its asymptotic expansion (in $t$) are always positive integers, nor that one should expect this. Galway \cite{Ga} proved this curious fact to be true    using alternating permutations and relations to Euler numbers. Stanley \cite{St}, answering a question of Galway, then gave a nice combinatorial interpretation of the coefficients in the asymptotic expansion (\ref{Ram}) as the number of fixed-point-free alternating involutions in the symmetric group $S_{2n}$, providing a second proof that the coefficients are positive integers.

Berndt and Kim \cite{BK} study more general partial theta functions.
\begin{align}\label{falseb}
f_b(t):=2\sum_{n=0}^\infty (-1)^n q^{n^2+bn},
\end{align} where  $b$ is real.
They show that, similarly to \eqref{Ram}, $f_b(t)$ admits an asymptotic expansion of the form
\begin{equation}\label{asexp}
f_b(t)\sim \sum_{n=0}^\infty a_n t^n,
\end{equation}
where the $a_n$ are given explicitely in term of Euler numbers and Hermite polynomials \cite[Theorem 1.1]{BK}.  For the purposes of this paper, we do not require the explicit shape of the $a_n$, so we do not state it here.  (We point out a small typo in \cite{BK} equation (1.3):  the exponent of $(1-t)/(1+t)$ should read $(1-2b)/4$ rather than $(2b-1)/4$, as the authors correctly state in \cite{BK} equation (2.9).)

In analogy to the results established by Galway and Stanley pertaining to the coefficients of (\ref{Ram}), Berndt and Kim prove that the coefficients $a_n$ of the generalized partial  theta functions defined in (\ref{falseb}) are integers if $b\in\mathbb N$ \cite[Theorem 2.5]{BK}, and make the following conjecture regarding their positivity.

\begin{conj}[Berndt-Kim \cite{BK}]\label{conjBerndtKim}
For any positive integer $b$, for sufficiently large $n$, the coefficients $a_n$ in the asymptotic expansion \eqref{asexp}, have the same sign.
\end{conj}
\noindent
The purpose of this note is to prove this conjecture of Berndt and Kim.

\begin{theorem}\label{maintheorem}
The Berndt-Kim conjecture is true.  More  precisely, if $b\equiv 1,2\pmod 4$, then $a_n>0$ for $n\gg0$, and if $b\equiv 0,3\pmod 4$, then $a_n<0$ for $n\gg 0$.
\end{theorem}

\section{Proof of Theorem \ref{maintheorem}}

\noindent
To prove the theorem, we first observe that  for integers $b\geq 2$,
\begin{equation}\label{recurse}
f_{b+1}(t)=-q^{-b}\Big(f_{b-1}(q)-2\Big).
\end{equation}
We  proceed by induction on $b$ to prove Theorem \ref{maintheorem}. The case $b=1$ follows from \cite{Ga} as mentioned in \S \ref{intro}. To prove the case $b=2$ we employ the fact that the $q$-series in this case is essentially a modular form.
 To be more precise, we have  that
\[
f_2(t)=-q^{-1}\left(g\left(\frac{i\theta}{2\pi}\right)-1\right),
\]
where we have adopted the notation $\theta = \log\left(\frac{1+t}{1-t}\right)$ from \cite{BK}.  The function $g\left(\tau\right)$, where $\tau \in \mathbb H$, is the modular form defined by
\[
g(\tau):=\sum_{n\in\Z}(-1)^n e^{2\pi i n^2\tau}=\frac{\eta(\tau)^2}{\eta(2\tau)},
\]
where $\eta(\tau):=e^{\frac{2\pi i \tau}{24}}\prod_{n\geq 1}(1-e^{2\pi i n \tau})$ is Dedekind's eta-function, a well known modular form of weight $1/2$. Employing the modular transformation of $\eta$ (see \cite{Radem} e.g.)
\[
\eta\left(-\frac{1}{\tau}\right)=\sqrt{-i\tau}\eta(\tau)
\]
we obtain that
\[
g\left(\frac{i\theta}{2\pi}\right)\to 0
\]
as $\theta\to 0^+$, and thus as $t\to 0^+$. Therefore
\[
f_2(t)\sim q^{-1} = \frac{1+t}{1-t},
\]
which clearly has positive coefficients in its $t$-expansion. \\
We now assume Theorem \ref{maintheorem} holds for some $b-1\geq 1$, and prove that it also holds for $b+1$.  (The first two inductive cases $b-1=1,2$ are proven above.)
We use \eqref{recurse} and split
\[
q^{-b}=(1+t)^{b}(1-t)^{-b}.
\]
Since the $t$-expansion of $(1+t)^{b}$ is finite and contains only positive  coefficients, it suffices to prove that the $t$-coefficents in the asymptotic expansion of
$$-(1-t)^{-b}\left(f_{b-1}(q)-2\right)$$ eventually all have the same sign.  We will address the fact that the sign is dependent on the residue class of $b\pmod{4}$ as stated in Theorem \ref{maintheorem} later.

It is easy to show that
\[
(1-t)^{-b}=\sum_{j=0}^\infty\binom{b-1+j}{b-1}t^j.
\]
By induction, we may assume without loss of generality that
\[
f_{b-1}(q)-2\sim\sum_{n\geq 0}\alpha_n t^n, \]
as $t\to 0^+$, where $\alpha_n>0$ for $n\geq L,$ for some $L\geq 0$. Set
\[
\beta_n:=\binom{b-1+n}{b-1}.
\]
It suffices to show that
\begin{align}\label{mixedsum}
\sum_{0\leq n\leq m}\alpha_n \beta_{m-n}
\end{align}
is positive for $m\gg 0$.  
We break the sum (\ref{mixedsum}) into two parts
\begin{align*}\Sigma_1 &:= \sum_{0\leq n \leq L-1}\alpha_n\beta_{m-n},
\ \ \ \ \ \ \ \ \Sigma_2:= \sum_{L\leq n \leq m}\alpha_n\beta_{m-n},
\end{align*} where we assume $m\gg 0$ is sufficiently large to ensure $m\geq L$.

Since $L$ is independent of $m$ and the $\beta_j$'s are monotonically decreasing (in $j$, for fixed $b$), we may bound
\begin{align*}
\Sigma_1 \ll
\beta_m\ll m^{b-1}, \end{align*}
as $m\to \infty$.
Similarly, $\Sigma_2$ may be estimated from below by
\[
\Sigma_2 \gg
\sum_{L\leq n\leq m}\beta_{m-n}=\sum_{0\leq n\leq m-L}\beta_n\gg \sum_{1\leq n\leq m}n^{b-1}
\gg\int_0^m x^{b-1}dx\gg m^{b},
\] as $m\to \infty$.  Since clearly  $m^{b}\gg m^{b-1},$ 
the positivity of $\Sigma_2$ dominaties, and we have that the coefficients in the asymptotic $t$-expansion of $f_{b+1}(q)$ have for sufficiently large $m$ the same sign as claimed. 

The more precise claim that the coefficients in the asymptotic $t$-expansion for $f_b(q)$ are eventually positive for $b\equiv 1,2\pmod{4}$ and eventually negative for $b\equiv 0,3\pmod{4}$ follows easily by induction using (\ref{recurse}), and the previously established facts that the coefficients in the asymptotic $t$-expansion of $f_1(q)$ and $f_2(q)$ are eventually positive.

\end{document}